# UNCOVERING A PALAEOTSUNAMI TRIGGERED BY MASS MOVEMENT IN AN ALPINE LAKE

M. Naveed Zafar[1,2], Denys Dutykh[3,4], Pierre Sabatier[2], Mathilde Banjan[2], and Jihwan Kim[5]

[1]Univ. Grenoble Alpes, Univ. Savoie Mont Blanc, CNRS, LAMA, Chambéry, France
[2]Univ. Savoie Mont Blanc, CNRS, EDYTEM, 73370 Le Bourget-du-Lac, France
[3]Mathematics Department, Khalifa University of Science and Technology, PO Box 127788, Abu Dhabi, United Arab Emirates
[4]Causal Dynamics Pty Ltd, Perth, Australia
[5]Instituto Português do Mar e da Atmosfera (IPMA), Lisbon, Portugal

## ABSTRACT

Mass movements and delta collapses are significant sources of tsunamis in lacustrine environments, impacting human societies enormously. Palaeotsunamis play an essential role in understanding historical events and their consequences along with their return periods. Here, we focus on a palaeo event that occurred during the Younger Dryas to Early Holocene climatic transition, *ca.,* 12,000 years ago in the Lake Aiguebelette (NW Alps, France). Based on high-resolution seismic and bathymetric surveys and sedimentological, geochemical, and magnetic analyses, a seismically induced large mass transport deposit with an initial volume of $767172\ m^3$ was identified, dated and mapped. To investigate whether this underwater mass transport produced a palaeotsunami in the Lake Aiguebelette, this research combines sedimentary records and numerical models. Numerical simulations of tsunamis are performed using a visco-plastic landslide model for tsunami source generation and two-dimensional depth-averaged nonlinear shallow water equations for tsunami wave propagation and inundation modelling. Our simulations conclude that this sublacustrine landslide produced a tsunami wave with a maximum amplitude of approximately $2\ m$ and run-up heights of up to $3.6\ m$. The modelled sediment thickness resulting from this mass transport corroborates well with the event deposits mapped in the lake. Based on our results, we suggest that this sublacustrine mass transport generated a significant tsunami wave that has not been reported previously to the best of our knowledge.

## 1. INTRODUCTION

The most common geophysical disturbances that trigger tsunamis are earthquakes, landslides, or their combinations. Tsunamis can vary widely in size, and their extent depends on the characteristics of these geophysical perturbations. Earthquake tsunamis are the most significant and global tsunamis, with large wavelengths of thousands of kilometres and relatively small amplitudes (Kanamori, 1972; Abe, 1973; Ward, 1980; Okal et al., 1988; Dutykh et al., 2006; Dutykh & Dias, 2007; Dias et al., 2014). They can devastate entire ocean regions with immense energy produced by both horizontal and vertical displacements of the bottom motion (Tanioka & Satake, 1996; Dutykh & Dias, 2009; Dutykh et al., 2012; Song et al., 2017), such as in two recent megathrust earthquake-tsunamis, the 2004 Indonesian (Lay et al., 2004) and the 2011 Tohoku (Simons et al., 2011) tsunamis. On the other hand, tsunamis caused by mass movement related to landslides (Hampton et al., 1996; Ward, 2001) or to volcanic collapse (Kienle et al., 1987; Pararas-Carayannis, 1992; Paris, 2015) are generally more localized; however, they have the potential to generate large tsunami waves (Bardet et al., 2003; Beizel et al., 2012; Dutykh & Kalisch, 2013; Løvholt et al., 2015; Waldmann, 2021). They are complicated and have a more diverse nature, although their wavelength is several hundred meters, much shorter than that of earthquake tsunamis (Harbitz et al., 2006; Dutykh et al., 2011b). Therefore, in comparison to earthquake tsunamis, prognostic analysis is less developed for mass-movement tsunamis, and thus, we have a limited understanding of this type of tsunami hazard. One of the main reasons for this lack of knowledge related to mass movement-induced tsunamis is that they fortunately have a very low recurrence rate, so there are limited records available (Blikra et al., 2005; Yamada et al., 2012; Kremer et al., 2015, 2012, 2021; Lane et al., 2017; Urgeles & Camerlenghi, 2013; Waldmann, 2021; Nigg, 2021).

Mass-movement tsunamis can be generated by subaerial or subaqueous mass transport deposits. Tsunamis triggered by subaqueous mass movement are receiving significant attention as natural hazards that pose a risk to coastal communities (Bardet et al., 2003; Satake, 2012). Strong regional coseismic displacements or volcanic eruptions can trigger subaqueous mass movements. Depending on the size and characteristics of the mass movement and the bathymetric setting, these mass failures may be tsunamigenic (Urgeles & Camerlenghi, 2013; Dutykh & Kalisch, 2013; Løvholt et al., 2015, 2017; Kim et al., 2019). One recent example of such a tsunami is the 1998 Papua New Guinea tsunami (Tappin et al., 2001; Synolakis et al., 2002; Lynett et al., 2003), generated by seismically triggered subaqueous mass movement. This catastrophic event, characterized by a maximum local runup of 15 meters, resulted in the tragic loss of over 2,000 lives. Other examples include the event of 1994 at Skagway, Alaska (Kulikov et al., 1996; Rabinovich et al., 1999; Thomson et al., 2001) and 1979 at Nice Airport, France (Assier-Rzadkieaicz et al., 2000; Dan et al., 2007). Due to their long recurrence intervals, studying palaeotsunamis is essential. Geological archives such as sediment deposits are the only solid evidence from which we can estimate both the frequency and magnitude of past tsunamis (Nanayama et al., 2003; Paris et al., 2007; Goto et al., 2011; Sugawara et al., 2014; Chagué-Goff et al., 2017; Biguenet et al., 2021; Cordrie et al., 2022).

The study of marine palaeotsunami sediments has been the focus of attention in recent decades, with numerous publications documenting the investigation of palaeotsunamis sediment archives at several marine sites. They investigated the recurrence rate (Monecke et al., 2008; Migeon et al., 2017; Kempf et al., 2017) and runup height (Bondevik et al., 1997; Petersone et al., 2011) and determined the inundation distance of past tsunami events using sedimentary records (Switzer et al., 2012; Abe, Goto, and Sugawara 2012; Goto et al., 2019; Ishimura & Yamada, 2019; Paris et al., 2023). However, investigating their lacustrine counterparts has received little attention thus far (Sammartini et al., 2019; Kremer et al., 2021; 2012). Mechanisms of tsunami generation by mass movements are not limited to the marine environment; they also occur in different lacustrine settings, *e.g.,* fjord-type lakes or fault-bounded basin lakes, which have large mass movements relative to the lake basin. Several historical lacustrine tsunamis events have been documented in various Alpine lakes in Switzerland, including tsunamis in Lake Geneva (Kremer et al., 2012), Lake Lauers (Bussmann & Anselmetti, 2010), Lake Lucerne (Hilbe & Anselmetti, 2014), Lake Brienz (Girardclos et al., 2007) and Lake Sils (Nigg et al., 2021). These earthquake-triggered lacustrine tsunami events were generated by delta collapses and subaqueous mass movements. Sediment archives are the only way to gain deeper insight into these events, and geoscientists can reconstruct these low-frequency events using sediment records. Several mass-movement events have been reconstructed in pre-Alpine lakes based on geophysical and sedimentological evidence that transport large amounts of sediments from lateral slopes to deep basins (Van Rensbergen et al., 1998; Chapron et al., 1999; Strasser & F. Anselmetti, 2008; Beck, 2009; Wilhelm et al., 2013). Some of these events occurred in the Younger Dryas–Early Holocene (YD-EH) transition time range, such as Lake Bourget and Lake Annecy (close to Lake Aiguebelette), and were studied using seismic reflection surveys. One of these types of mass-movement events was discovered in Lake Aiguebelette in the YD-EH climate transition with an initial mass of 767172 $m^3$ (Banjan et al., 2022). This event was identified thanks to high-resolution seismic and bathymetric data associated with sedimentary, geochemical, and magnetic proxies and has a seismic origin enhanced by rapid climate change and glacial retreat, such as crustal rebound and erosional unloading at this climate transition (Banjan et al., 2022).

Numerical modelling of tsunamis serves as an alternative method for estimating potential tsunami hazards when direct measurements are absent (Dutykh et al., 2011a; LeVeque et al., 2011; Kim et al., 2019; Cordrie et al., 2022). The numerical modelling of palaeotsunamis is an essential part of multidisciplinary research that bridges the gap between the practical applications of tsunami sediments and geological studies (Jaffe et al., 2016). This research aims to conduct the first numerical study of the French-Alpine Lake Aiguebelette and to answer whether this reconstructed subaqueous mass movement could generate tsunamis. From the sedimentary characterization of this sublacustrine landslide as an input parameter, we reproduce the tsunami wave by solving the nonlinear shallow water equations using finite volume schemes and computing the tsunami wave generation, propagation, and inundation. This research study is also relevant for other peri-Alpine lakes that exhibit similar sublacustrine mass movements (Rapuc et al., 2022) and associated tsunami hazards.

The present manuscript is organized as follows. A description of the study site and sedimentary sequence is provided in Section 2. Section 3 introduces numerical models of mass movement and tsunamis. Simulation results for mass transport deposits, tsunami propagation and inundation are presented in Section 4. We discuss our results, particularly the validation of the simulations with geological observations in Section 5. The paper concludes in the final Section 6, summarizing our current research and outlining future perspectives.

## 2. STUDY SITE AND LAKE SEDIMENTARY SEQUENCES

Lake Aiguebelette is a pre-Alpine lake that formed after glaciers retreat. It is in the northwestern Alps, France, at an elevation of 373 $m$ above sea level with a surface area of 5.45 $km^2$. The lake's average depth is 30 $m$, and the maximum depth of the deepest basin is approximately 70 $m$. The 'Leysse de Novalaise' is the main tributary entering the lake's northern part. The geology of the catchment area (see Fig. 1A) and an aerial picture of the lake are shown in Fig. 1B.

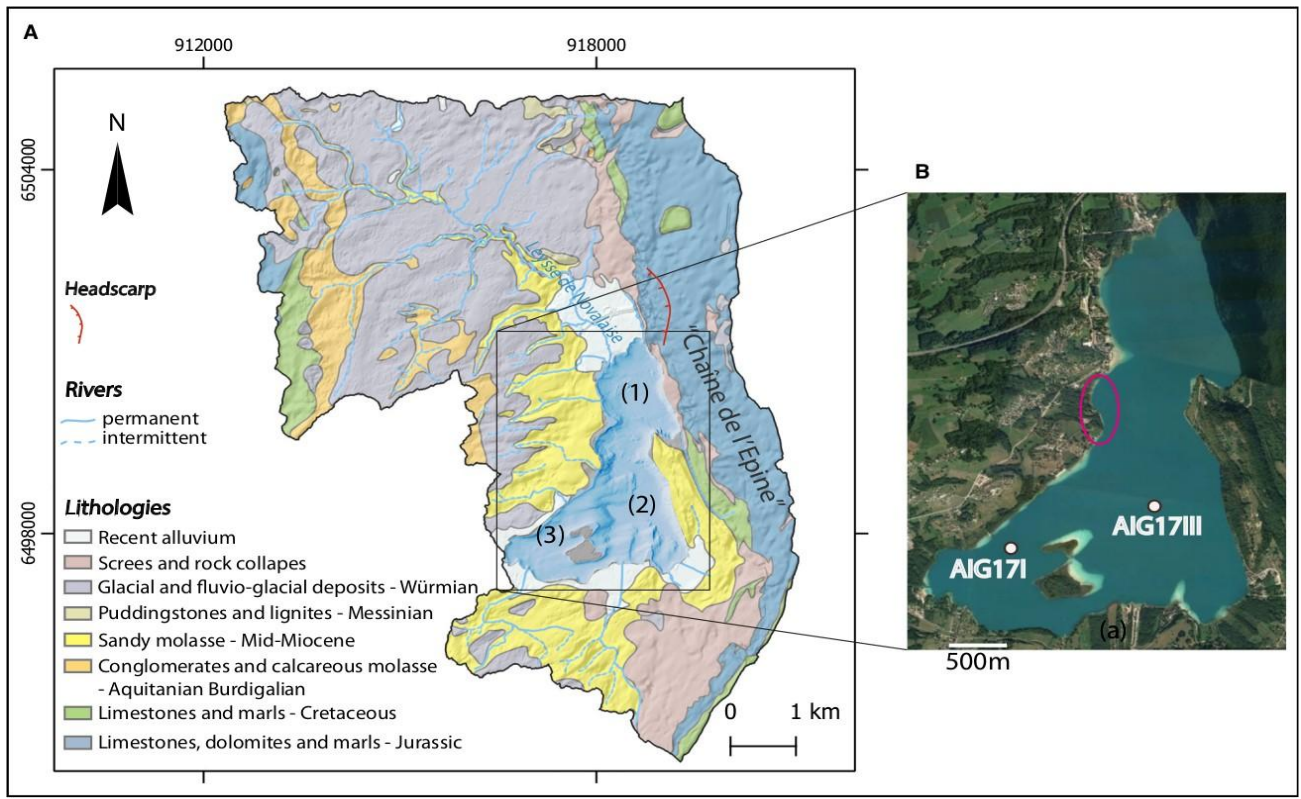

**Figure 1:** (A) Geological map of the Lake Aiguebelette catchment with the three subbasins of the lake: (i) northern basin, (ii) deep basin and (iii) shallow southern basin. (B) Depicts the two core sites of the AIG17I and AIG17III sediment sequences in the shallow and deep basins, respectively, published in Banjan et al., 2022. The red circle on the western shore of the lake shows the sliding mass position. We can observe the littoral carbonate platform on the western coast in light blue, except for the area of sliding mass in the red circle.

Banjan et al., 2022, present a high-resolution bathymetric survey and seismic data acquisition, which allows for identifying a large mass wasted deposit. The 16 m long main core named AIG17III belongs to the deep basin of Lake Aiguebelette ($45.5503°N$, $5.8015°E$) at a depth of $70\ m$, allowing the characterization of a 115 cm thick turbidite–homogenite dated at the Younger Dryas-Early Holocene climatic transition. A seismic survey enabled to determine the thickness of the final mass-transport deposits (ranging from $0.6\ m$ to $1.7\ m$) in areas that are not gas blind (highlighted in Fig. 2A), where the channelized mass transport deposits indicate that mass deposits come from an area on the western flank of the lake (Fig. 2B) where carbonates are absent; however, littoral carbonates are seen everywhere on the western side of the lake, except where the channelized morpho-bathymetric feature arises (Fig. 1B). The initial sliding mass was therefore identified on the western shore of the lake by mapping sediment head traces through GIS analysis of bathymetric data (Fig. 2B). By reconstructing a prefailure slope profile on the bathymetric data, the volume of this mass transport deposit was estimated to be approximately $767172\ m^3$ using kriging interpolation (the area used for kriging and the shape of the initial mass with thickness are shown in Fig. 2B).

## 3. NUMERICAL MODELS USED IN THIS STUDY

### 3.1 Mass-movement Model:

Sublacustrine mass movements involve complex layered structures, including a dense bottom debris layer covered by a thin layer entrained by turbidity currents. This bottom layer can be taken as either a clay-rich visco-plastic flow or a granular, incoherent flow. There will be no deformation in the case of viscoplastic flow unless the shear stress exceeds the yield strength. When it exceeds this limit, the flow begins to behave as a shear-thinning, non-Newtonian fluid. This type of flow behaviour is taken as Herschel-Bulkeley rheology, which is incorporated in an applied numerical depth-averaged mass transport model called BingClaw (Løvholt et al., 2017; Kim et al., 2019) modified from (Huang & García, 1997; Imran et al., 2001). The Herschel-Bulkley rheological model for simple shear can be defined as follows:

$$|\dot\gamma/\dot\gamma_r|^n = \begin{cases} sgn(\dot\gamma)\left((\tau/\tau_y) - 1\right) & if \ |\tau| > \tau_y, \\ 0 & Otherwise, \end{cases}$$

where $\tau$ and $\tau_y$ are the shear stress and yield strength, respectively, $\dot\gamma$ is the strain rate and $\dot\gamma_r$ is the reference strain rate, given as

$$\dot\gamma_r = \left(\frac{\tau_y}{\mu}\right)^{\frac{1}{n}}.$$

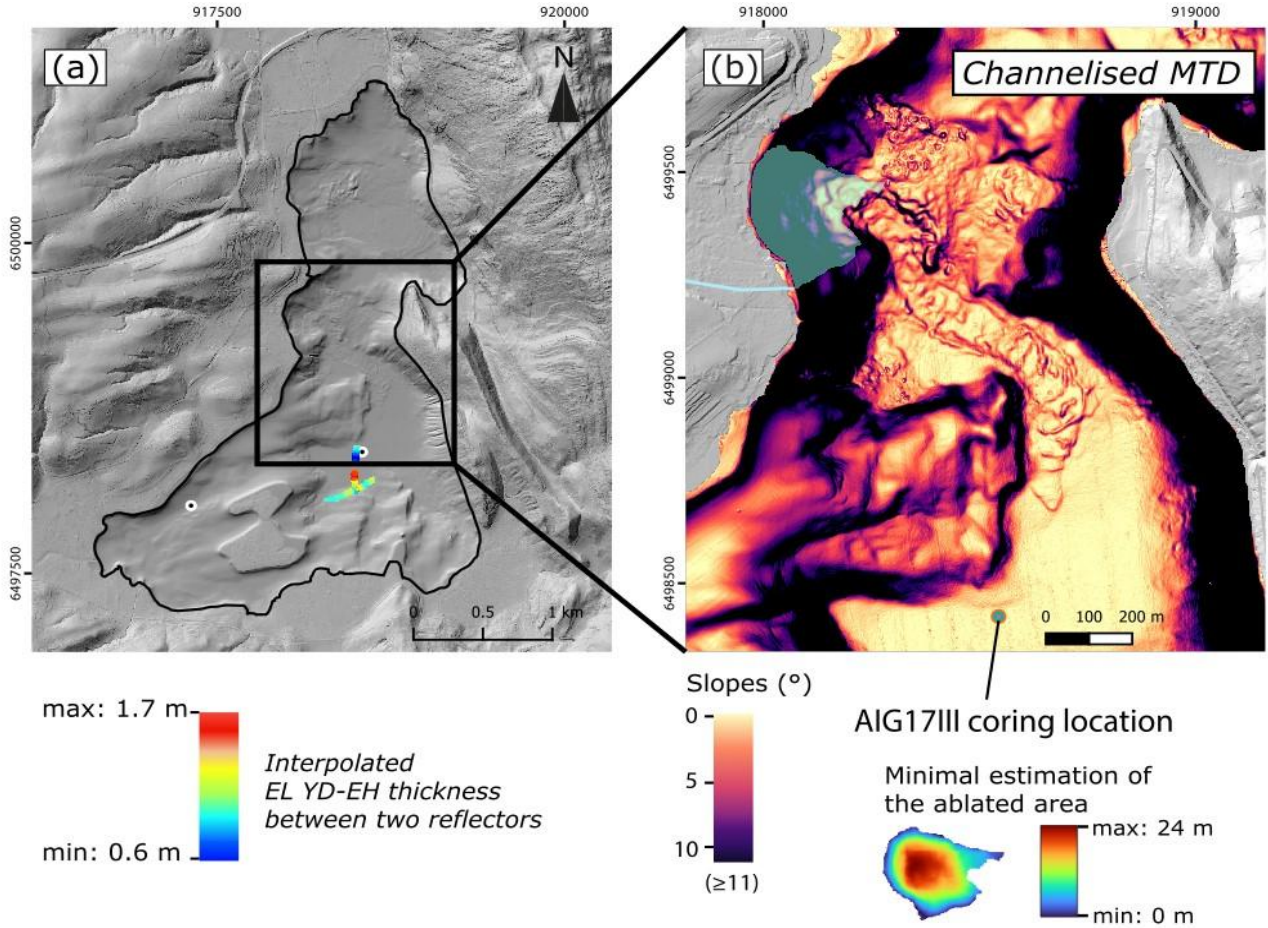

**Figure 2:** (A) Thickness of the final mass-transport deposits (ranging from 0.6 m to 1.7 m) mapped from the seismic profile. (B) A map of channelized mass transport deposits and slope angles is presented, and the initial sliding mass is highlighted on the western flank of the lake. The thickness of the initial mass (ranging from 0 to 24 m) is shown at the bottom of the figure.

The model integrates the mass balance equations over the entire flow depth and solves the momentum equations separately for the top plug layer, which exhibits no shear deformation, and the bottom shear layer. The entire flow depth is shown as $H = H_s + H_p$; here, $H_p$ and $H_s$ represent the thickness of the plug and the shear layers in the vertical direction, respectively. Volumetric flux through the plug layer ($H_p \boldsymbol{u}_p$) and the shear layer ($H_s \boldsymbol{u}_s$) with velocities $\boldsymbol{u}_p$ and $\boldsymbol{u}_s$ parallel to the slope within the plug and the shear layer, respectively:

$$\frac{\partial}{\partial t}(H_p + H_s) + \nabla \cdot (H_p \boldsymbol{u}_p + H_s \boldsymbol{u}_s) = 0,$$

$$\left(1 + Cm\frac{\rho_w}{\rho_d}\right)\frac{\partial}{\partial t}(H_p \boldsymbol{u}_p) + \nabla \cdot (H_p \boldsymbol{u}_p \boldsymbol{u}_p)$$

$$+ g' H_p \nabla (H_p + H_s + b) + \boldsymbol{u}_p \left(\frac{\partial}{\partial t} H_s + \nabla \cdot (H_s \boldsymbol{u}_s)\right) = -\frac{\boldsymbol{u}_p}{||\boldsymbol{u}_p||}\frac{\tau_y + \tau_d}{\rho_d},$$

$$\left(1 + Cm\frac{\rho_w}{\rho_d}\right)\frac{\partial}{\partial t}(H_s \boldsymbol{u}_s) + \nabla \cdot (\alpha \ H_s \boldsymbol{u}_s \boldsymbol{u}_s)$$

$$+ g'H_s \nabla(H_p + H_s + b) - \boldsymbol{u}_p \left(\frac{\partial}{\partial t} H_s + \nabla \cdot (H_s \boldsymbol{u}_s)\right) = -\frac{\boldsymbol{u}_p}{||\boldsymbol{u}_p||} \frac{\tau_y f_s}{\rho_d}.$$

This depth-averaged model operates in 2HD (or two horizontal dimensions), and all vector quantities and the nabla operator $\nabla$ are two-dimensional horizontal directions. Here, $b(x, y)$ represents the bathymetric depth, $C_m$ is the added-mass coefficient, $\rho_w$ and $\rho_d$ are the density of ambient water and the landslide, respectively, $\alpha$ is the velocity form factor in the shear layer, and t is the time coordinate. The reduced gravitational acceleration (considering the buoyancy effects) is defined as $g' = g(1 - \rho_w - \rho_d)$, and $\tau_d$ is the viscous drag at the free surface, split into friction drag $\tau_f$ and a pressure drag term $\tau_p$ given as:

$$\tau_f = \tfrac{1}{2} C_F \rho_w \boldsymbol{u}_p ||\boldsymbol{u}_p||, \qquad \tau_p = \tfrac{1}{2} C_P \rho_w max(0, -\boldsymbol{u}_p \cdot \nabla H)\boldsymbol{u}_p,$$

where $C_F$ and $C_P$ represent the hydrodynamic coefficients of friction and pressure drags, respectively. At the right-hand side of the share-layer momentum balance, $\tau_y f_s$ represents the net shear stress at the bed; here,

$$f_s = \beta \cdot \left(\frac{||\boldsymbol{u}_p||}{\gamma_r H_s}\right)^n,$$

with shape factor $\beta = \left(1 + \frac{1}{n}\right)^n$.

In BingClaw, dominant terms are modelled using finite volume and source terms by finite difference methods. This framework utilizes the conservation law package ClawPack (Mandli et al., 2016). In computational grids where the earth pressure does not exceed the material shear strength, the cell remains static. If it does, the dynamic equations are solved using the Godunov fractional step method. ClawPack first solves the equations without friction terms using the finite volume method and then accounts for frictional terms in the next fractional step. For a more detailed derivation, see (Løvholt et al., 2017; Kim et al., 2019).

### 3.2 Tsunami Model:

We use the GeoClaw model, a version of the Clawpack that is a set of numerical methods for solving nonlinear hyperbolic conservation laws. GeoClaw solves the two-dimensional depth-averaged nonlinear shallow water equations using well-balanced Godunov-type finite volume schemes on rectangular grids. It solves the Riemann problem at each cell interface and allows the discontinuous solution or shock wave that arises in hyperbolic PDEs. Additionally, we can incorporate adaptive mesh refinement for robust and efficient solutions for large-scale geophysical problems, and it captures wet–dry interfaces efficiently for inundation modelling. For a detailed description of GeoClaw, see (LeVeque et al., 2011). We used the temporal volumetric moving bathymetry resulting from the mass-movement model BingClaw as input for tsunami generation and computed the tsunami propagation and inundation for Lake Aiguebelette. The results are presented in the next section.

## 4. RESULTS

### 4.1 Mass-movement Simulations

The volume and shape of the displaced mass are estimated by kriging interpolation using the current bathymetry due to insufficient information on the prefailure bathymetry. However, with the help of geological and seismic data, the comparison of the simulation results with the observations provides the ground for the validation of the results and leads us to a realistic volume. Mass movement simulations are conducted on a $(480 \times 400)$ rectangular computational grid with a cell size of approximately 5 meters and by using high-resolution bathymetry data with a $2\ m$ resolution. The range of the non-Newtonian fluidity index $n$ is $0 < n < 1$ for shear thinning Herschel-Bulkley flow, and we used $n = 1/2$ for the clay-silty sediments in the simulations. Other key parameters, such as landslide density, yield strength, added mass and initial volume, are listed in Table 1.

Table 1: Model parameters for the mass-movement simulations.

| Volume | Rheological model | Slide Density | Fluid Density | Initial Yield Strength |
|---|---|---|---|---|
| $767172\ m^3$ | Herschel-Bulkley | $1500\ kg/m^3$ | $1000\ kg/m^3$ | $65\ pa$ |

Typical values for the pressure force $C_P = 1.0$ and the frictional drag $C_F = 0.01$ were also considered for the simulation of mass failure to maintain realistic velocity and acceleration. The thickness and run-out distance of mass-movement deposits at various times are shown in Fig. 3. Mass deposits take approximately 12 minutes to travel to the final deposit location.

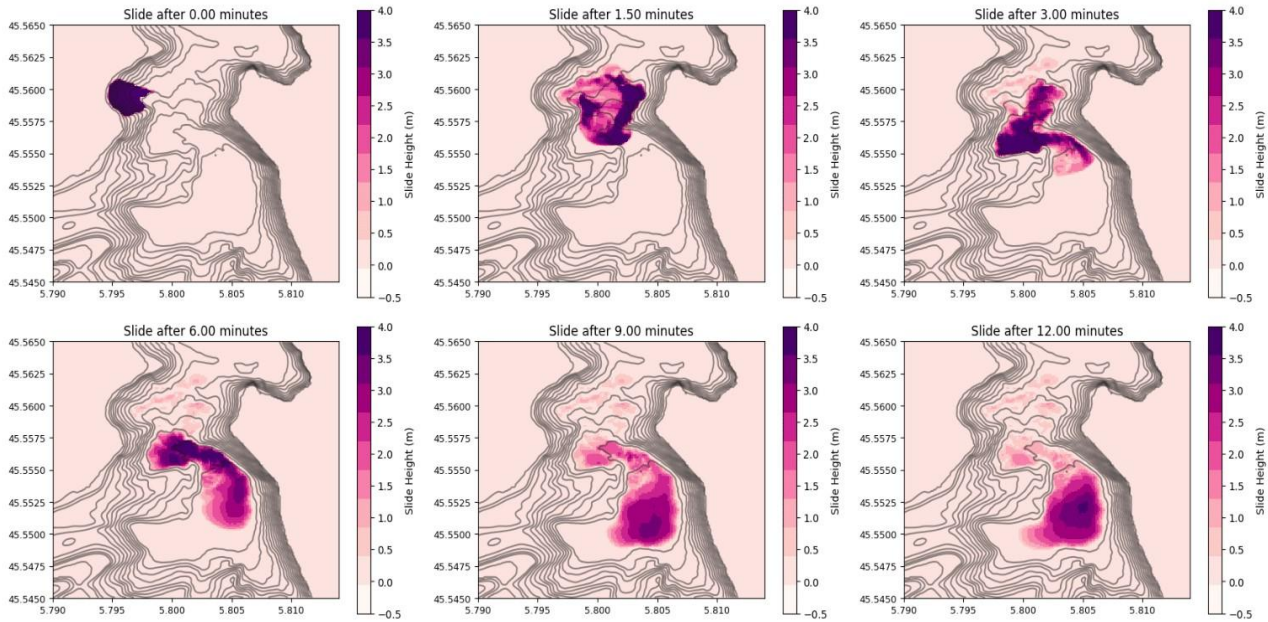

**Figure 3:** Snapshots of the thickness and runout distance of the mass-movement numerical model results at 0, 1.5, 3, 6, 9, and 12 min.

### 4.2 Tsunami Simulations

In this study, the primary source of tsunami generation is the temporal variations in lakebed bathymetry caused by sublacustrine mass movement. Such progressive changes in bathymetry are the outputs of the BingClaw mass movement model and are included as inputs to the tsunami model. These perturbations are directly transformed into the overlying water column and the water surface. Tsunami propagation and inundation are performed using the nonlinear shallow water equations (NSWEs) for Lake Aiguebelette using an offshore bathymetry dataset with a resolution of 1.2 meters and onshore high-resolution topography less than 1 meters employed for inundation modelling. The rectangular computational grids ($720 \times 600$) were used with a cell size of approximately 5 meters.

#### *4.2.1 Wave Propagation and inundation:*

The propagation of tsunami waves over time as predicted by NSWE is shown in Fig. 4. Black dots represent offshore and onshore synthetic gauge numbers from $1 - 10$ and $11 - 15$, respectively. Upon failure of the mass, the surface amplitude is rapidly displaced by approximately 5 $m$ initially at $t = 6\ s$. After that, the leading wave propagated ~800 $m$ from the impact area in the direction of the sliding mass with a maximum amplitude greater than 1.5 m at t = 30 s and reached the opposite eastern shore, right in front of the massive failure. After 30 $sec$, the majority of the wave spread to the northeastern, while a lesser portion heads towards the southeastern. The wave propagates across the entire lake within 2 minutes, although approximately 12 minutes are needed for the mass deposits to reach the final deposition area. We can observe that this sublacustrine mass movement produced a significant tsunami wave with immense wavelength.

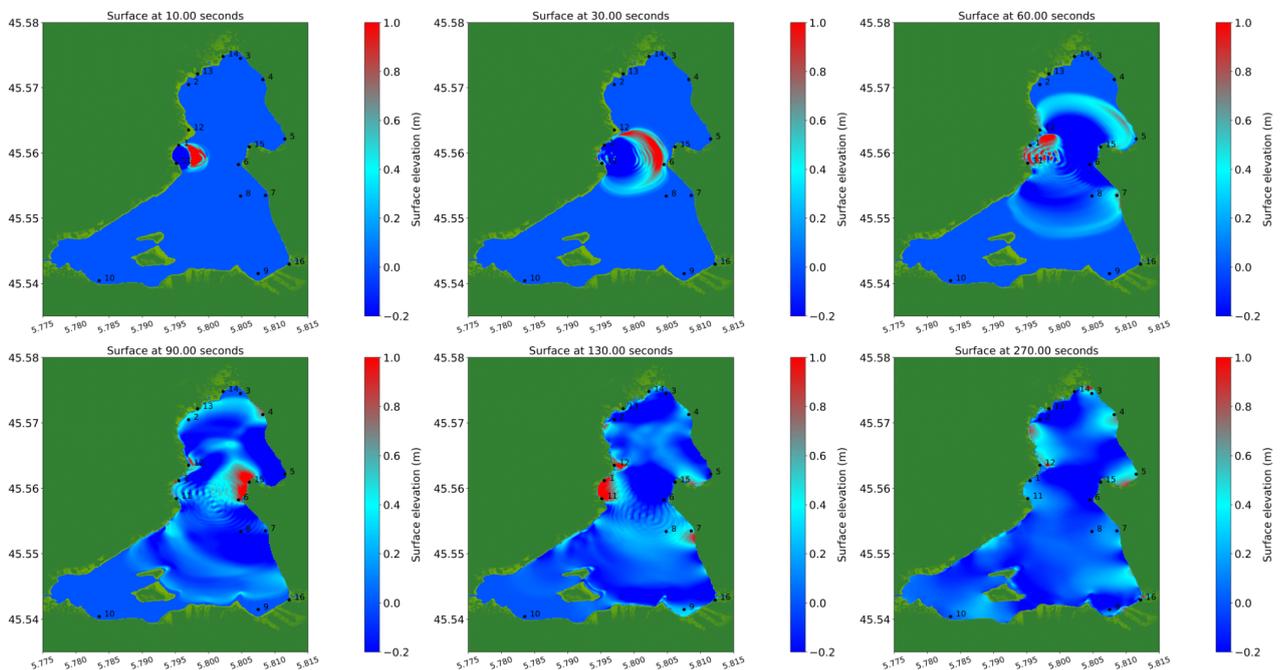

**Figure 4:** Snapshots of palaeotsunami wave propagation at 10, 30, 60, 90, 130, and 270 sec. The locations of synthetic gauges 1-15 are indicated with black dots.

Time series of water surface elevation at synthetic gauges were extracted from the tsunami simulations and are shown in the figure together with the maximum amplitude across the entire lake (Fig. 5). Waves of approximately 2 meters amplitude are observed near the source on gauge 1 and 1.5 meters near the opposing shore on gauge 6 in the direction of the leading wave.

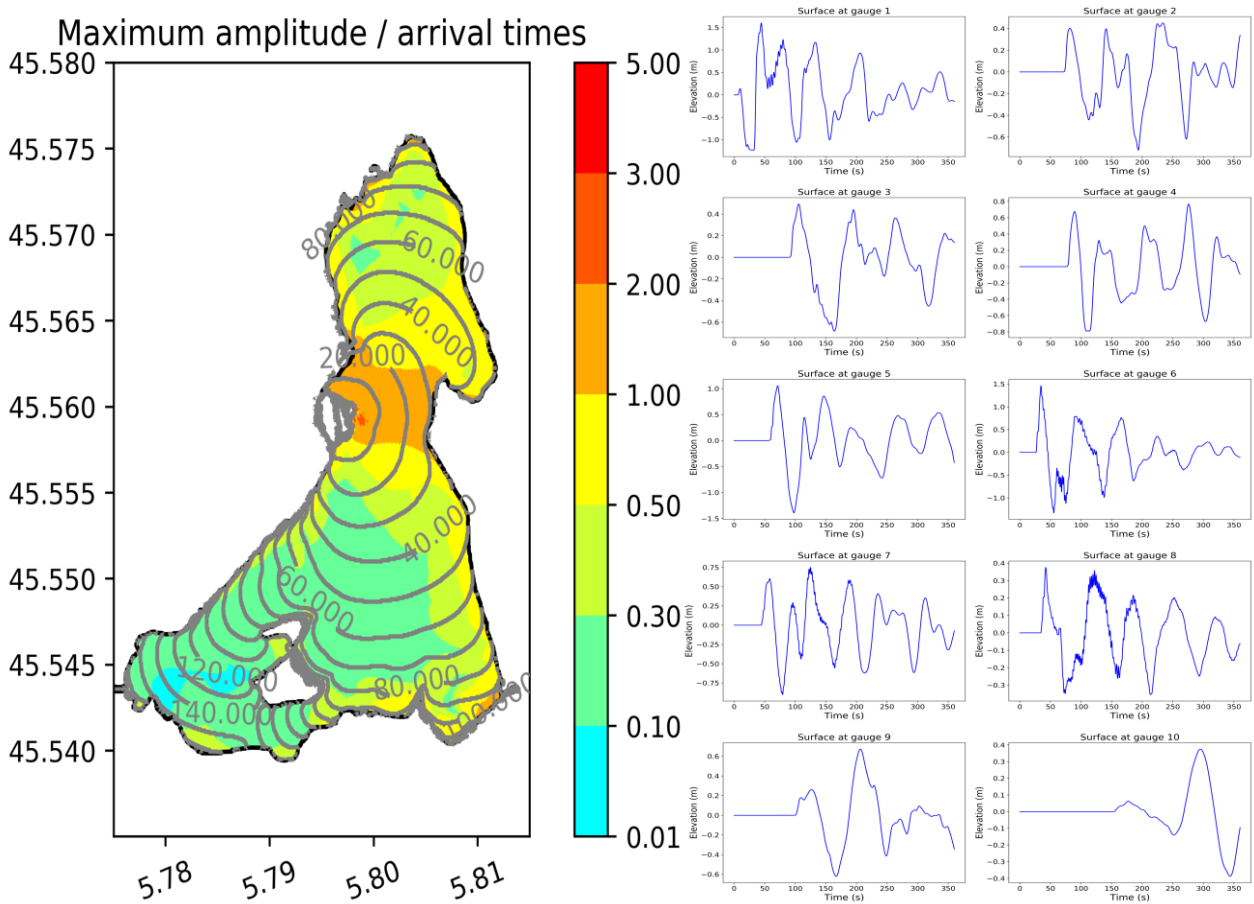

**Figure 5:** Maximum wave amplitude across the lake over time (left) with time series of wave amplitude at synthetic gauges 1-10 at offshore locations (right).

At the lake deepest basin, the amplitude is approximately 0.5 meters on gauge 8. In the figure of maximum amplitude over arrival times (Fig. 5), we can observe maximum amplitudes throughout the lake, with higher amplitudes in the centre of the lake in the east–west direction. We also conducted inundation modelling at the finer grids by using the current topography data of the lake shore to observe the effects of tsunami waves at nearby shores. The tsunami runup with synthetic gauges at onshore locations can be observed in Fig. 6. The topography surrounding the lake features significantly steeper slopes, particularly on the east side of the lake, which is a mountain area, which is why we do not observe the significant runout distance by using the current topography data. The west shore has a maximum tsunami runup distance of approximately 150 meters and a maximum tsunami runup height of 3.6 meters.

## 5. DISCUSSION

Simulating sublacustrine mass movement numerically provides insight into past events, allows the reconstruction of palaeo events and identifies future risks. Simulations of such landslides are complex because of their dynamic nature. Tsunamis generated by underwater landslides highly depend on the characteristics of the initial sliding mass and bathymetry (Harbitz et al., 2006; Lovholt et al., 2015; Watts, 1998). In this study, the post-bathymetry data serve as our foundation. We acknowledge that substantial bathymetric changes may have occurred over the last 12,000 years, introducing uncertainties to the initial sliding mass and bathymetry. Good observation data can constrain the uncertainties; bathymetry surveys and sediment analyses of the Lake Aiguebelette provide detailed information about the initial failure mass and physical parameters such as density, water content, and sediment grain sizes about palaeo-sublacustrine mass movement. Multiple tests are conducted to determine the remaining essential parameters, including the initial yield strength, crucial to determining the runout distance (Kim et al., 2019). Seismic and sedimentological data enabled us to validate our mass movement model and reduced the error between the mass-movement simulations and the geological observations. Based on the parameters given in Table (1), our numerical simulations align well with the observation of the Younger Dryas – Early Holocene event. An excellent agreement has been found, and we obtained the same sediment thickness and runout distance derived from seismic data and simulation at the final deposit location, as illustrated in the Fig. 7. Due to the remarkable agreement between the computed results and the observation, we are confident that the mass movement simulations can be used to identify the potential tsunami at the Lake Aiguebelette.

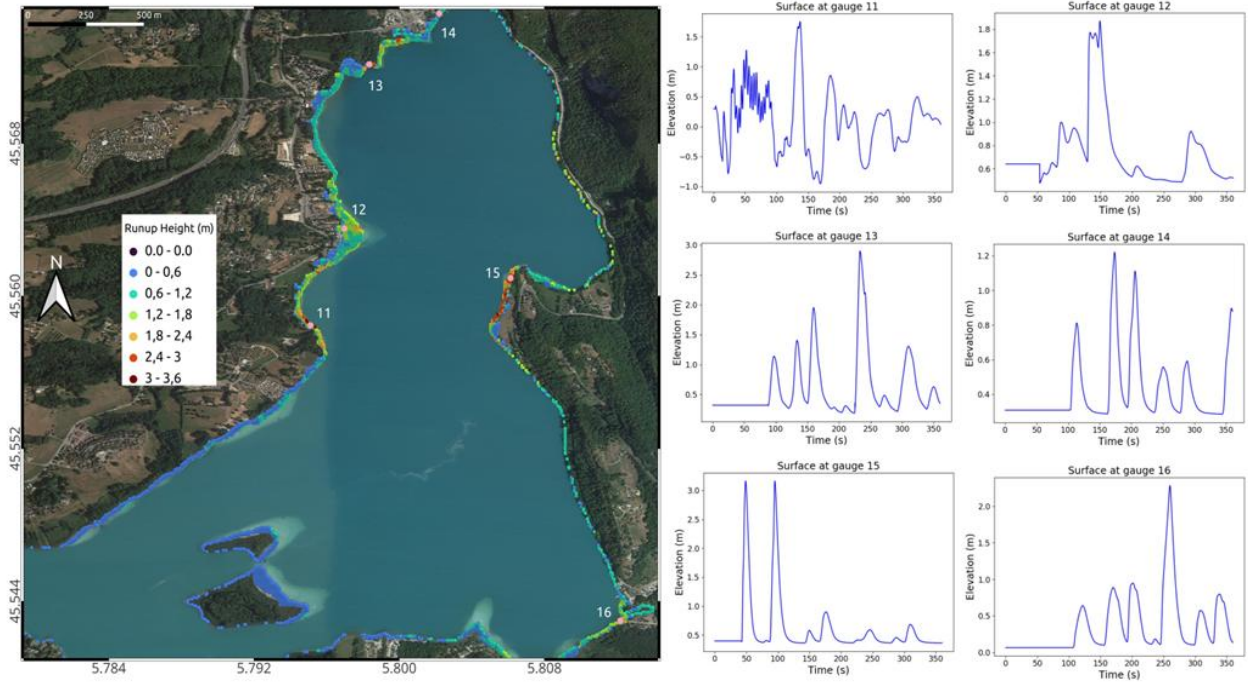

**Figure 6:** Wave height (left) and time series of runup heights at the indicated synthetic gauges 11-16 (right) across the lake shore, bearing in mind the present-day bathymetry and shoreline.

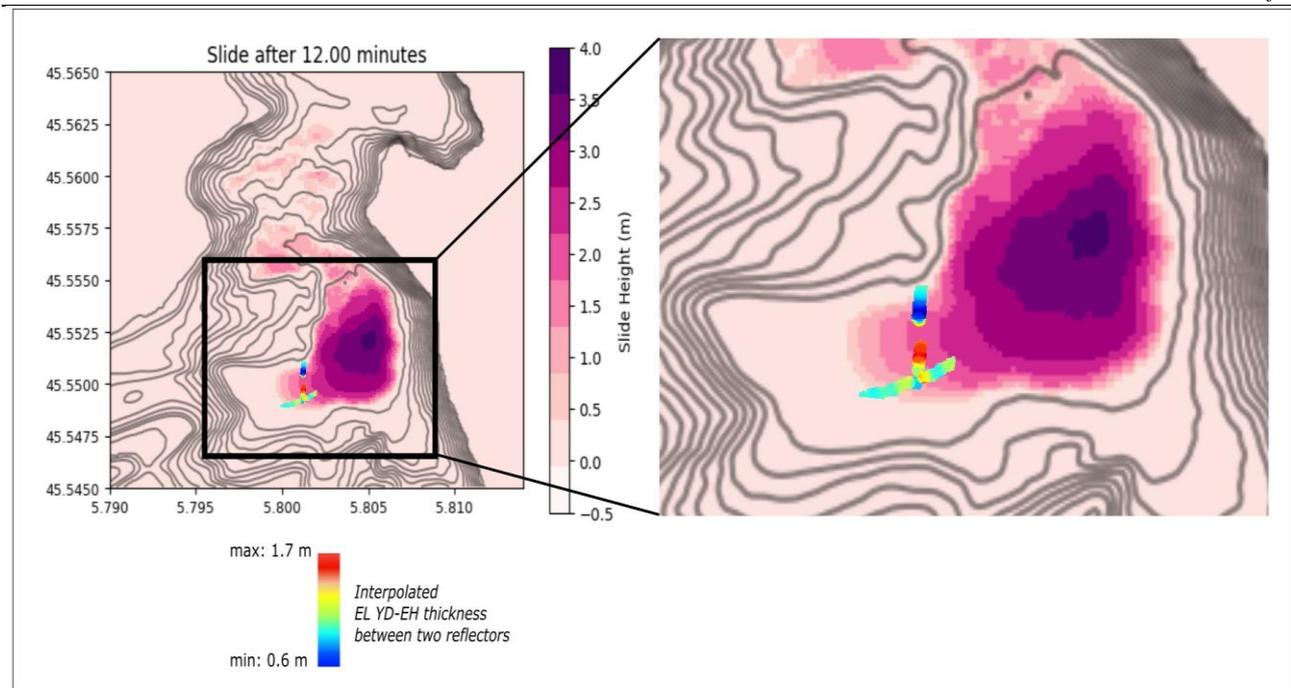

**Figure. 7:** Validation of mass deposit thickness with seismic profile at the final depositional location.

Using the well-constrained mass movement model, our investigation showed that seismically induced mass movement resulted in significant palaeotsunami waves in the Lake Aiguebelette. We can observe from the mass-movement propagation (Fig. 3) and wave propagation (Fig. 4) that the direction of the leading tsunami wave toward the east follows the direction of the mass deposits in the early stages of the collapse. However, over time, the speed of the mass deposits diminishes, and its influence on the tsunami wave's direction becomes less pronounced in the deceleration phase. It is interesting to observe that all the mass transport moves toward the south after 1.5 minutes, but the majority of the leading wave moves toward the north. It seems that tsunami wave direction is determined by the initial phase of mass-movement deposits and the bathymetry features of the area. Since tsunami wave direction is essential for locating potential tsunami deposits, the results of this numerical study provide insight into potential palaeotsunami deposit locations. Eventually, we intend to locate the tsunami deposits in the direction of the highest positive wave to gain a deep understanding of this palaeo-sublacustrine tsunami event. It is crucial to note that we simulate the tsunami event that occurred approximately 12,000 years ago, even though water levels fluctuated over time, while this variation in water levels might influence our findings. Regarding tsunami hazards in this area, we performed the calculations for the run-up distance and height by using the current topography, as shown in the Fig. 6. The most affected area is along the western side of the lake because it is less steep than the eastern shore, although the run-up height is greater at the eastern shore. Run-up is mainly dependent on the local bathymetric features close to the shoreline (*e.g.,* Grilli et al., 2009; De Blasio 2011). Overall, we do not obtain a significant runup distance because the steep slope confines the inundation process.

## 6. CONCLUSION AND PERSPECTIVES

In this paper, we modelled the mass movement and a related tsunami in Lake Aiguebelette using a viscoplastic Herschel-Bulkley rheological model coupled in one way with the hydrodynamic nonlinear shallow water equation model. This research study provides evidence for the potential of significant sublacustrine mass movement-induced palaeotsunami in an alpine lake, and our numerical simulations shed light on how palaeosublacustrine mass movement in a lake environment can lead to tsunamis. Based on detailed geologic and bathymetric data, numerical simulations of mass movement match well with observed deposits. The simulated runout distance and thickness are consistent with the geological observations. For the tsunami simulations, we observe that the sublacustrine mass movement with an initial volume of 767172 $m^3$ could have transported enough sediments to generate a basin-wide tsunami in Lake Aiguebelette with a 3.6 meter runup height. Tsunami simulations show that the average amplitude of the tsunami wave in the east–west direction is approximately 2 meters, which causes a substantial tsunami. This multidisciplinary research is crucial to understanding the past event that occurred in the Younger Dryas–Early Holocene climatic transition and if such event occurred in the future climate change how could be impact the shore of the lake.

Our future research will focus on analysing the dispersive effects of uncovered palaeotsunamis in Lake Aiguebelette. We will conduct a comparative study between the nonlinear shallow water equation and the dispersion model to observe the importance of dispersion effects in small-scale lake environments. Additionally, we are in the process of organizing a survey to locate onshore tsunami deposits around the Lake Aiguebelette. The identification and analysis of these deposits will allow us to deepen our understanding of palaeotsunami studies and will contribute to reducing the uncertainties associated with palaeotsunami studies.

## ACKNOWLEDGEMENTS

The acquisition of bathymetric data would not have been possible without Conservatoire des Espaces Naturels de Savoie, CCLA and Reserve Naturelle Regionale du Lac d'Aiguebelette. This publication is based upon work supported by the Khalifa University of Science and Technology under Award No. FSU-2023-014. MNZ gratefully acknowledges the financial support received from the MESRI and ED MSTII in the form of a PhD scholarship to conduct research for his doctoral thesis.

## REFERENCES AND CITATIONS